\documentclass[10pt]{scrartcl}
\usepackage{enumitem}
\usepackage{amsmath}
\usepackage{amssymb}
\usepackage{mathtools}

\usepackage{amsthm}
\usepackage{abstract}
\usepackage{xcolor}
\usepackage{comment}
\usepackage{dsfont}
\usepackage{hyperref}
\usepackage{amsmath}
\DeclareMathOperator*{\argmax}{arg\,max}

\usepackage{algorithm}
\usepackage{algpseudocode}
\hypersetup{
    colorlinks=true,
    linkcolor=blue,
    filecolor=magenta,      
    urlcolor=cyan
    }
\makeatletter
\newcommand*{\barfix}[2][.175ex]{%
  \mathpalette{\@barfix{#1}}{#2}%
}
\newcommand*{\@barfix}[3]{%
  \vbox{%
    \kern#1\relax
    \hbox{$#2#3\m@th$}%
  }%
}
\makeatother

\newtheorem{theorem}{Theorem}
\newtheorem{thm}{Theorem}[section]
\newtheorem{corollary}[thm]{Corollary}
\newtheorem{lemma}[thm]{Lemma}

\newcommand{\footremember}[2]{%
    \footnote{#2}
    \newcounter{#1}
    \setcounter{#1}{\value{footnote}}%
}

\usepackage[margin=1in]{geometry}
\title{\vspace{-2cm}Heavy and Light Paths and Hamilton Cycles}  
\author{%
Sahar Diskin \footremember{alley}{School of Mathematical Sciences, Tel Aviv University, Tel Aviv 6997801, Israel. Email: sahardiskin@mail.tau.ac.il.}%
\and Dor Elboim \footremember{trailer}{Princeton University, Princeton NJ, United States. Email: delboim@princeton.edu.}%
}
\begin{document}
\maketitle
\begin{abstract}
    Given a graph $G$, we denote by $f(G,u_0,k)$ the number of paths of length $k$ in $G$ starting from $u_0$. In graphs of maximum degree 3, with edge weights $i.i.d.$ with $exp(1)$, we provide a simple proof showing that (under the assumption that $f(G,u_0,k)=\omega(1)$) the expected weight of the heaviest path of length $k$ in $G$ starting from $u_0$ is at least 
    \begin{align*}
        (1-o(1))\left(k+\frac{\log_2\left(f(G,u_0,k)\right)}{2}\right),
    \end{align*}
    and the expected weight of the lightest path of length $k$ in $G$ starting from $u_0$ is at most
    \begin{align*}
        (1+o(1))\left(k-\frac{\log_2\left(f(G,u_0,k)\right)}{2}\right).
    \end{align*}

    We demonstrate the immediate implication of this result for Hamilton paths and Hamilton cycles in random cubic graphs, where we show that typically there \textit{exist} paths and cycles of such weight as well. Finally, we discuss the connection of this result to the question of a longest cycle in the giant component of supercritical $G(n,p)$.
\end{abstract}
\section{Introduction and main results}
\textit{Stochastic programming} refers to a general class of optimisation problems where uncertainty is modelled by a probability distribution on the input variables. In one family of such problems, one considers a graph $G=(V,E)$ with a random weight function $w: E(G)\to \mathbb{R}$, and asks for the minimum or maximum weight object. There has been extensive research into these types of questions. To name a few, Janson \cite{J99} studied the weighted distance in the complete graph $K_n$ with edge weights $i.i.d.$ with $exp(1)$; Aldous \cite{A01} considered the expectation of the minimal total weight of a perfect matching in $K_{n,n}$ with exponential weights with rate 1; Beveridge, Frieze, and McDiarmid \cite{BFM98} studied the length of the minimum spanning tree in regular graphs with random weights; and, even the question of a long cycle in the giant component of supercritical $G(n,p)$ can be seen as a question of finding heavy paths and cycles in graphs with random weights (see Section \ref{discussion} for further discussion and references on the matter). We refer the reader to \cite{UP13} for a collection of several papers and results in stochastic optimisation, and to Chapter 20 of \cite{FK16} for many more optimisation problems in graphs with random weights.

In this note, we propose and analyse a method to study heavy and light paths (and cycles) of given length, starting from a given vertex, in graphs with random weights. Given a graph $G=(V,E)$, we denote by $f(G,u_0,k)$ the number of paths of length $k$ in $G$, starting from $u_0$. Before stating our main result, let us heuristically explain its relation to the function $f(G,u_0,k)$. Indeed, given a graph $G$ with random weight function, if $u_0$ has only one path of length $k$ starting from it, we do not anticipate to be able to obtain paths whose weight deviates from the expected weight of a path of length $k$. On the other hand, if there are many paths of length $k$ starting from $u_0$, one might hope to use that in order to find paths whose weight does deviate from the expected weight. Intuitively, the more paths of length $k$ there are, one might hope for a possibly larger deviation. 

The main result of this note is as follows:
\begin{theorem}\label{expectation}
Let $G$ be a graph with $\Delta(G)\le 3$. Let $k$ be an integer. Let $u_0\in V(G)$ with $f(G,u_0,k)=\omega(1)$. Let $w:E(G)\to \mathbb{R}$ be a weight function on the edges of $G$, assigning independently to each edge a weight distributed according to $exp(1)$. Let $H(u_0,k)$ be the weight of a heaviest path of length $k$ in $G$ starting from $u_0$, and let $L(u_0, k)$ be the weight of a lightest path of length $k$ in $G$ starting from $u_0$. Then,
\begin{align*}
    \mathbb{E}[H(u_0,k)]\ge \left(1-o(1)\right)\left(k+\frac{\log_2\left(f(G,u_0,k)\right)}{2}\right),
\end{align*}
and
\begin{align*}
    \mathbb{E}[L(u_0,k)]\le \left(1+o(1)\right)\left(k-\frac{\log_2\left(f(G,u_0,k)\right)}{2}\right).
\end{align*}
\end{theorem}
Let us try to give some intuition behind this result. Consider the following process, starting from a vertex $u_0$ in a graph $G$ of maximum degree $3$. At each step, we observe the remaining edges of the last vertex in the path we constructed thus far (there are at most two such edges, since we assumed our graph has maximum degree $3$). If we have two choices allowing us to complete a path of length $k$, we can choose the maximum between them, or the minimum between them, having an expected weight $\frac{3}{2}$ and $\frac{1}{2}$, respectively. Otherwise, we must choose the only edge possible allowing us to complete a path of length $k$, having an expected weight $1$. Heuristically, we anticipate the choice between two edges to be possible at $\log_2\left(f(G,u_0,k)\right)$ of the vertices. We can therefore hope to construct in this (greedy) manner a heavy path of weight $$\frac{3}{2}\cdot\log_2\left(f(G,u_0,k)\right)+\left(k-\log_2\left(f(G,u_0,k)\right)\right)=k+\frac{\log_2\left(f(G,u_0,k)\right)}{2},$$
and a light path of weight
$$\frac{1}{2}\cdot\log_2\left(f(G,u_0,k)\right)+\left(k-\log_2\left(f(G,u_0,k)\right)\right)=k-\frac{\log_2\left(f(G,u_0,k)\right)}{2}.$$
Thus, in a sense, Theorem \ref{expectation} validates the above heuristics.

As an example demonstrating these heuristics and the result, consider the complete binary tree of depth $n$, $T_n$, with root $r$. We have that $f\left(T_n, r, n\right)=2^n$, as, indeed, at every step we have two choices, and we thus may anticipate to gain $\frac{n}{2}$ in the expected weight. And indeed, by Theorem \ref{expectation}, we have that the expected weight of a heaviest path of length $n$ in $T_n$ starting from the root $r$ is at least $\left(1-o(1)\right)\frac{3n}{2}$, and the expected weight of a lightest path of length $n$ in $T_n$ starting from the root $r$ is at most $\left(1+o(1)\right)\frac{n}{2}$.

Furthermore, given a mild requirement on the length of the path, we are able to show the typical \textit{existence} of such paths in these graphs.
\begin{theorem} \label{general statement}
Let $G$ be a graph on $n$ vertices with $\Delta(G)\le 3$. Let $\omega\left(n^{\frac{1}{2}}\right)\le k \le n$ be an integer. Let $u_0\in V(G)$ with $f(G,u_0,k)=\omega_n(1)$. Assign independently to each edge of $G$ a weight distributed according to $exp(1)$. Then, \textbf{whp}\footnote{With high probability, that is, with probability tending to $1$ as $n$ tends to $\infty$.}, the following holds.
\begin{itemize}
    \item[(1)]\label{heavy path} There exists a path of length $k$ in $G$ starting from $u_0$ with weight at least $$\left(1-o(1)\right)\left(k+\frac{\log_2\left(f(G,u_0,k)\right)}{2}\right);$$
    \item[(2)]\label{light path} There exists a path of length $k$ in $G$ starting from $u_0$ with weight at most $$\left(1+o(1)\right)\left(k-\frac{\log_2\left(f(G,u_0,k)\right)}{2}\right).$$
\end{itemize}
\end{theorem}

It is worth noting that Theorem \ref{general statement} utilises two greedy algorithms, whose (quite simple) analysis shows that they typically construct paths of these weights.

One immediate implication concerns Hamilton paths in random cubic graphs. Indeed, when $G$ is a random 3-regular graph on $n$ vertices, it is known that \textbf{whp} the graph is Hamiltonian (\cite{RW92}). In fact (see, \cite{RW92, FJMRW96}), denoting by $X$ the random variable counting the number of Hamilton cycles in $G$, we have that \textbf{whp} $\frac{1}{n}\cdot \left(\frac{2}{\sqrt{3}}\right)^n\le X \le n\left(\frac{2}{\sqrt{3}}\right)^n.$ Therefore, together with Theorem \ref{general statement}, we immediately obtain the following corollary:
\begin{corollary}
Let $G\sim G_{n,3}$ be a random $3$-regular graph. Assign independently to each edge of $G$ a weight distributed according to $exp(1)$. Then, \textbf{whp} there exists a Hamilton path in $G$ of weight at least $\left(1-o(1)\right)\left(1+\frac{\log_2\left(\frac{2}{\sqrt{3}}\right)}{2}\right)n>1.1n.$ Furthermore, \textbf{whp} there exists a Hamilton path in $G$ of weight at most $\left(1+o(1)\right)\left(1-\frac{\log_2\left(\frac{2}{\sqrt{3}}\right)}{2}\right)n<0.89n.$
\end{corollary}
Using the above, we will be able to construct Hamilton \textit{cycles} with such weights in random cubic graphs:
\begin{corollary} \label{cubic graphs}
\sloppy Let $G\sim G_{n,3}$ be a random $3$-regular graph. Assign independently to each edge of $G$ a weight distributed according to $exp(1)$. Then, \textbf{whp} there is an Hamilton cycle in $G$ of weight at least $\left(1-o(1)\right)\left(1+\frac{\log_2\left(\frac{2}{\sqrt{3}}\right)}{2}\right)n$. Furthermore, \textbf{whp} there is an Hamilton cycle in $G$ of weight at most $\left(1+o(1)\right)\left(1-\frac{\log_2\left(\frac{2}{\sqrt{3}}\right)}{2}\right)n$.
\end{corollary}

A few comments about our notation are in place. When considering a path $P$, we denote by $|P|$ the number of vertices in the path. Similarly, when we mention the length of the path, we mean the number of vertices in it. Given a graph $G$ of maximum degree $\Delta(G)$ and a vertex $u_0\in V(G)$, we assume its neighbours are numbered $v_1,\cdots ,v_\ell$, where $0 \le \ell \le \Delta(G)$. Furthermore, given an integer $k$, we denote by $\alpha_i$ (where $i\in[\ell]$) the fraction of paths of length $k$ in $G$ starting from $u_0$ and going through the edge $u_0v_i$ (note that $\sum_{i=1}^{\ell}\alpha_i=1$). That is, $f(G\setminus u_0, v_i, k-1)=\alpha_i\cdot f(G, u_0, k)$. We denote by $u_P$ the last vertex of a path $P$. Finally, we denote by $P\cup\{v_i\}$ the concatenation of the edge $(u_Pv_i)$ to $P$. We omit rounding signs for the clarity of presentation.

The structure of the note is as follows. The proofs throughout the note utilise two greedy algorithms, which are presented In Section \ref{Greedy algorithms}. In Section \ref{key proof}, we prove Theorems \ref{expectation} and \ref{general statement}. In Section \ref{cubic} we prove Corollary \ref{cubic graphs}. Finally, in Section \ref{discussion} we discuss the relation of our results to the question of estimating a longest cycle in supercritical $G(n,p)$ and possible variants of the algorithms.

\section{The greedy algorithms}\label{Greedy algorithms}
We begin with the description of Heavy Greedy Algorithm (HGA). At every step the algorithm chooses the edge that maximises $w(v_Pv_i)+\frac{\log_2\alpha_i}{2}$, concatenates it to the path $P$, and continues with $G\setminus v_p$. Formally:
\begin{algorithm}[H]
\caption{Heavy Greedy Algorithm}\label{alg:HGA}
    \hspace*{\algorithmicindent} \textbf{Input}: A graph $G$, a weight function $w:E(G)\to\mathbb{R}_{\ge0}$, $u_0\in V(G)$ with $f(G, u_0,k)\ge 1$, and an integer $k$.
    \begin{algorithmic}[1]
    \State $P_0=\{u_0\}$, $G_0=G$, $m=0$.
    \While {$|P_m|<k$}
        \State $j=\argmax_{i}\left\{w(u_{P_m}v_i)+\frac{\log_2\alpha_i}{2}\right\}$. \label{alg: heavy choose}
        \State $P_{m+1}=P_m\cup\{v_j\}$.
        \State $G_{m+1}=G_m\setminus u_{P_m}$.
        \State $m\to m+1$.
    \EndWhile
    \State $\tilde{P}=P_m$.
    \end{algorithmic}
    \hspace*{\algorithmicindent} \textbf{Output}: A path $\tilde{P}$.
\end{algorithm}

The description of Light Greedy Algorithm (LGA) is very similar. At step 3, instead of choosing the edge that maximises $w(v_Pv_i)+\frac{\log_2\alpha_i}{2}$, it simply chooses the edge that minimises $w(v_Pv_i)-\frac{\log_2\alpha_i}{2}$. We thus omit the full description.

We note that in both algorithms, at every step $m$, the $\alpha_i$ are with respect to the graph $G_m$ and to the vertex $u_{P_m}$.

Observe that both algorithms will always output $\tilde{P}$ with $|\tilde{P}|=k$. Indeed, an edge which will not allow it to continue to a path of length $k$ will have $\alpha_i=0$, and hence $\log_2\alpha_i\to-\infty$, and therefore will not be chosen neither in HGA nor in LGA. In particular, at every step $m$ (before the algorithm stops), we have that $d_{G_m}(v_{P_m})\ge 1$.

We denote by $HGA(G, w, u_0, k)$ and $LGA(G, w, u_0, k)$ the path obtained by HGA and LGA, respectively, with these parameters. 

Finally, we note that these algorithms take as a blackbox the values of $f(H,u,k')$ for various subgraphs $H\subseteq G$, vertices $u\in V(G)$ and values of $k'$. While for some graphs, mainly graphs rising from random processes, it is possible to estimate the value of these functions, this task may prove difficult in the case of general graphs. Still, there are several techniques which one may employ in order to estimate this function through sampling (see, for example, Chapter 11 of \cite{MU17}).

\section{Proof of Theorems \ref{expectation} and \ref{general statement}}\label{key proof}
In order to further simplify the analysis, if we receive $G$ and $u_0$ such that $d_{G}(u_0)=3$, we first remove the edge with smallest $\alpha_i$ amongst the edges incident to $u_0$, and otherwise we do not change the graph. We denote the resulting graph by $G'$. We then have that $f(G',u_0,k)\ge f(G,u_0,k)/3$. We run HGA and LGA on $G'$ instead of on $G$, and as we will see, this will have negligible effect on the weight of the paths which are constructed.

We begin with the following lemma, bounding the expected weight of the path constructed by the algorithm. Given a graph $G$ with a weight function $w: E(G)\to\mathbb{R}_{\ge0}$ and a path $P$ in $G$, we denote by $w(P)$ the sum of the weight of the edges of $P$.

\begin{lemma}\label{l: expectation}
Let $G$ be a graph on $n$ vertices with $\Delta(G)\le 3$. Let $1\le k \le n$ be an integer. Let $w:E(G)\to \mathbb{R}$ be a weight function on the edges of $G$, assigning independently to each edge a weight distributed according to $\exp(1)$. 
\begin{itemize}
    \item[(1)]\label{l:heavy expectation} Let $P_H=HGA(G,w,u_0,k)$. Then, $\mathbb{E}[w(P_H)]\ge k+\frac{\log_2\left(f(G',u_0,k)\right)}{2}-1$;
    \item[(2)]\label{l:light expectation} Let $P_L=LGA(G,w,u_0,k)$. Then, $\mathbb{E}[w(P_L)]\le k-\frac{\log_2\left(f(G',u_0,k)\right)}{2}$.
\end{itemize}
\end{lemma}
\begin{proof}
We begin with the first item. We prove by induction on $k$, for all possible choices of graphs with maximum degree $3$. When $k=1$, we have only one path (of length $1$) of weight $0$, and we have that $f(G',u_0,1)=1$. Indeed, $0\ge 1+\frac{0}{2}-2$, as required.

Assume the hypothesis holds for $k-1$, and consider $k\ge 2$. If $d_{G'}(u_0)=1$, we have that HGA adds the edge $u_0v_1$ to the path. Noting that $\Delta(G'\setminus u_0)\le \Delta(G)\le 3$, we obtain by the induction hypothesis that
\begin{align*}
    \mathbb{E}[w(P_H)]&\ge \mathbb{E}\left[w(u_0v_1)+w\left(HGA\left(G'\setminus u_0, w, v_1, k-1\right)\right)\right]\\
        &=\mathbb{E}[w(u_0v_1)]+\mathbb{E}\left[w\left(HGA\left(G'\setminus u_0, w, v_1, k-1\right)\right)\right]\\
        &\ge 1+\left((k-1)+\frac{\log_2\left(f(G\setminus u_0, v_1, k-1)\right)}{2}-1\right)\\
        &=k+\frac{\log_2\left(f(G\setminus u_0, v_1)\right)}{2}-1, 
\end{align*}
where the last equality follows since $f(G', u_0, k)=\alpha_1\cdot f(G'\setminus u_0, v_1, k-1)$ and $\alpha_1=1$.
Otherwise, by construction, $u_0$ has degree two, and note that $\alpha_2=1-\alpha_1$. Let us denote by $\mathcal{A}$ the event that $w(u_0v_1)+\frac{\log_2\alpha_1}{2}\ge w(u_0v_1)+\frac{\log_2(1-\alpha_1)}{2}$. Then, by the definition of HGA, we have that:
\begin{align*}
    w(P_H)\ge &\max\left\{w(u_0v_1)+\frac{\log_2\alpha_1}{2}, w(v_0w_2)+\frac{\log_2(1-\alpha_1)}{2}\right\}\\
    &+\mathds{1}_{\mathcal{A}}\cdot\left(-\frac{\log_2\alpha}{2}+w\left(HGA\left(G'\setminus u_0, w, v_1, k-1\right)\right)\right)\\
    &+\mathds{1}_{\mathcal{A^C}}\cdot\left(-\frac{\log_2(1-\alpha)}{2}+w\left(HGA\left(G'\setminus u_0, w, v_2, k-1\right)\right)\right).
\end{align*}
Let us set $g(\alpha_1)\coloneqq \mathbb{E}\left[\max\left\{w(u_0v_1)+\frac{\log_2\alpha_1}{2},w(u_0v_2)+\frac{\log_2(1-\alpha_1)}{2}\right\}\right].$
Then, by the induction hypothesis, we obtain that
\begin{align*}
    \mathbb{E}[w(P_H)]\ge &g(\alpha_1)+\mathbb{P}(\mathcal{A})\cdot\left(-\frac{\log_2\alpha_1}{2}+(k-1)+\frac{\log_2\left(\alpha_1\cdot f(G',u_0,k)\right)}{2}-1\right)\\
    &+\mathbb{P}(\mathcal{A}^C)\cdot\left(-\frac{\log_2(1-\alpha)}{2}+(k-1)+\frac{\log_2\left((1-\alpha_1)f(G',u_0,k)\right)}{2}-1\right)\\
    &=g(\alpha_1)+\left(\mathbb{P}(\mathcal{A})+\mathbb{P}(\mathcal{A}^C)\right)\left(k-2+\frac{\log_2\left(f(G',u_0,k)\right)}{2}\right)\\
    &=g(\alpha_1)-1+k+\frac{\log_2\left(f(G',u_0,k)\right)}{2}-1.
\end{align*}
Thus, in order to complete the induction step, it suffices to show that $g(\alpha_1)-1\ge 0$. We may assume WLOG that $\alpha_1\le \frac{1}{2}$. We thus obtain that
\begin{align*}
    g(\alpha_1)=\frac{\log_2(\alpha_1)}{2}+\mathbb{E}\left[\max\left\{w(u_0v_1),w(u_0v_2)+\frac{\log_2(1-\alpha_1)-\log_2\alpha_1}{2}\right\}\right].
\end{align*}
Let us set $c(\alpha_1)\coloneqq \frac{\log_2(1-\alpha_1)-\log_2\alpha_1}{2}$. 

Now, we are interested in the expectation of $Z=\max\left\{X_1, X_2+c(\alpha_1)\right\}$ where $X_1,X_2\sim exp(1)$. We have that
\begin{align*}
    \mathbb{E}[Z]&=\int_0^{\infty}\mathbb{P}(Z>z)dz=\int_0^{\infty}\left(1-\mathbb{P}(Z\le z)\right)dz\\
    &=c(\alpha_1)+\int_{c(\alpha_1)}^{\infty}\left(1-\mathbb{P}(Z\le z)\right)dz=c(\alpha_1)+\int_{c(\alpha_1)}^\infty\left(1-\mathbb{P}(X_1\le z)\mathbb{P}(X_2\le z-c(\alpha_1))\right)dz\\
    &=c(\alpha_1)+\int_{c(\alpha_1)}^{\infty}\left(1-\left(1-\exp(-z)\right)\left(1-\exp\left(-x+c(\alpha_1)\right)\right)\right)dz\\
    &=c(\alpha_1)+\int_{c(\alpha_1)}^{\infty}\left(\exp(-z)+\exp(-z+c(\alpha_1))-\exp(-2z+c(\alpha_1))\right)dz\\
    &=c(\alpha_1)+1+\frac{\exp(-c(\alpha_1))}{2}.
\end{align*}
Therefore, $g(\alpha_1)=\frac{\log_2(1-\alpha_1)}{2}+1+\frac{\left(\frac{\alpha_1}{1-\alpha_1}\right)^{\frac{1}{\ln(4)}}}{2}$. 

We claim that $g(\alpha )\ge 1$ for all $0\le \alpha \le 1/2$, and hence $\mathbb{E}[w(p_H)]\ge k+\frac{\log_2\left(f(G',u_0,k)\right)}{2}-1,$ as required. To see this, note that $g(0)=g(1/2)=1$ and that $g(1/4)>1$. Thus, it suffices to show that there is a unique $0<\alpha <1/2$ such that $g'(\alpha )=0$. To this end, we solve 
\begin{equation}
    0=g'(\alpha )=\frac{2\alpha -(\alpha /(1-\alpha ))^{1/\ln 4}}{\alpha (\alpha -1) \ln 16}
\end{equation}
and therefore $2^{\ln 4} \alpha ^{\ln 4-1}=1/(1-\alpha )$. The last equation has at most two solutions in $(0,1)$ as the first function is concave and the second is convex. Since $1/2$ is a solution it means that there is at most one additional solution in $(0,1/2)$ but there has to be exactly one since $g(0)=g(1/2)$. This shows that $g(\alpha )\ge 1$ for all $0\le \alpha \le 1/2$ as needed.

We now turn to the second item, which follows an almost identical analysis to the first part. We will thus allow ourselves to be brief. Once again, we prove by induction on $k$ for all possible choices of graphs with maximum degree $3$. If $k=1$, we have one path of weight $0$, and $0\le 1-\frac{0}{2}$. 

Assume the hypothesis holds for $k-1$, and consider $k\ge 2$. If $d_{G'}(u_0)=1$, we have that LGA adds the edge $u_0v_1$ to the path. Noting that $\Delta(G'\setminus u_0)\le \Delta(G)\le 3$, we obtain by the induction hypothesis and calculations similar to before that
\begin{align*}
    \mathbb{E}[w(P_L)]&\le \mathbb{E}\left[w(u_0v_1)+w\left(LGA\left(G'\setminus u_0, w, v_1, k-1\right)\right)\right]\\
        &=k+\frac{\log_2\left(f(G\setminus u_0, v_1)\right)}{2}. 
\end{align*}
Otherwise, by construction, $u_0$ has degree two. We now set $\mathcal{B}$ to be the event that $w(u_0v_1)-\frac{\log_2\alpha_1}{2}\le w(u_0v_1)-\frac{\log_2(1-\alpha_1)}{2}$. Then, by the definition of LGA, we have that:
\begin{align*}
    w(P_L)\le &\min\left\{w(u_0v_1)-\frac{\log_2\alpha_1}{2}, w(v_0w_2)-\frac{\log_2(1-\alpha_1)}{2}\right\}
   \\& +\mathds{1}_{\mathcal{B}}\cdot\left(\frac{\log_2\alpha}{2}+w\left(LGA\left(G'\setminus u_0, w, v_1, k-1\right)\right)\right)\\&+\mathds{1}_{\mathcal{B^C}}\cdot\left(\frac{\log_2(1-\alpha)}{2}+w\left(LGA\left(G'\setminus u_0, w, v_2, k-1\right)\right)\right).
\end{align*}
Let us set $$\tilde{g}(\alpha_1)\coloneqq \mathbb{E}\left[\min\left\{w(u_0v_1)-\frac{\log_2\alpha_1}{2},w(u_0v_2)-\frac{\log_2(1-\alpha_1)}{2}\right\}\right].$$
Then, by the induction hypothesis and by calculations similar to before, we have that
\begin{align*}
    \mathbb{E}[w(P_L)]\le\tilde{g}(\alpha_1)-1+k-\frac{\log_2\left(f(G',u_0,k)\right)}{2}.
\end{align*}
Thus, in order to complete the induction step, it suffices to show that $\tilde{g}(\alpha_1)-1\le 0$. We are then interested in the expectation of $\tilde{Z}=\min\left\{X_1-\frac{\log_2\alpha}{2}, X_2-\frac{\log_2(1-\alpha_1)}{2}\right\}$ where $X_1,X_2\sim exp(1)$. We may assume WLOG that $\alpha_1\le \frac{1}{2}$. We then have that
\begin{align*}
    \mathbb{E}\big[\tilde{Z} \big]&=\int_0^{\infty}\mathbb{P}(\tilde{Z}>z)dz=\int_0^{\infty}\mathbb{P}\left(X_1\ge z+\frac{\log_2\alpha_1}{2}\right)\mathbb{P}\left(X_2\ge z+\frac{\log_2(1-\alpha_1)}{2}\right)dz\\
    &=0+\int_{-\frac{\log_2\alpha_1}{2}}^{\infty}\exp\left(-2z-\frac{\log_2\alpha_1+\log_2(1-\alpha_1)}{2}\right)\\
    &=\frac{\exp\left(\log_2\alpha_1-\frac{\log_2\alpha_1+\log_2(1-\alpha_1)}{2}\right)}{2}=\frac{\exp\left(-c(\alpha_1)\right)}{2}.
\end{align*}
Therefore, $\tilde{g}(\alpha_1)=\frac{\left(\frac{\alpha_1}{1-\alpha_1}\right)^{\frac{1}{\ln(4)}}}{2}.$ The derivative of $\tilde{g}(\alpha_1)$ is non-negative in the interval $\left[0,\frac{1}{2}\right]$. Thus the function achieves its maximum in this interval when $\alpha_1=\frac{1}{2}$, and we then have that $\tilde{g}(\alpha_1)<1$, as required.
\end{proof}

Note that Theorem \ref{expectation} follows immediately from the above Lemma. Indeed, we assume that $f(G,u_0,k)=\omega(1)$, and hence $k\pm\frac{\log_2\left(f(G',u_0,k)\right)}{2}=\left(1+o(1)\right)\left(k\pm\frac{\log_2\left(f(G',u_0,k)\right)}{2}\right)$. With this lemma at hand, we are also ready to prove Theorem \ref{general statement}.
\begin{proof}[Proof of Theorem \ref{general statement}]
We start by proving the first part of the theorem. Note that any graph $G$ on $n$ vertices with maximum degree $3$ has at most $\frac{3n}{2}$ edges. Enumerate the edges by $e_1,\dots ,e_m$ where $m\le \frac{3n}{2}$ and let $w_1,\dots ,w_m$ be the corresponding weights on these edges. Let $X=f(w_1,\dots ,w_m)$ be the weight of the heaviest path starting from $u_0$ of length $k$. For any $j\le m$ we define $X_j :=f(w_1,\dots ,w_{j-1},w_j',w_j,\dots ,w_m)$, where $w_j'$ is an $\exp(1)$ random variable independent of $w_1,\dots ,w_m$. That is, $X_j$ is defined in the same manner as $X$ when resampling only the weight on the edge $e_j$. By the Efron-Stein inequality \cite{steele1986efron} we have
\begin{equation}
    \text{Var}(X)\le \sum _{j=1}^m \mathbb E \big[ (X-X_j)^2 \big] \le \sum _{j=1}^m \mathbb E \big[ (w_j-w_j')^2 \big] \le 4m\le 6n,
\end{equation}
where in the second inequality we used that the weight of the heaviest path cannot change by more than $|w_j-w_j'|$ when changing $w_j$ to $w_j'$. By the first part of Lemma \ref{l: expectation} we have that $\mathbb{E}[X]\ge k+\frac{\log_2\left(f(G,u_0,k)\right)}{2}-1$ and therefore by Chebyshev's inequality for all $\epsilon >0$ we have 
\begin{equation}
    \mathbb P \Big( X \le k+\frac{\log_2\left(f(G,u_0,k)\right)}{2} -\epsilon k\Big) \le \mathbb P \big( |X-\mathbb E[X]|\ge \epsilon k-1  \big) \le \frac{6n}{(\epsilon k-1)^2}  \to 0,\quad n\to \infty,
\end{equation}
where we used the fact that $k=\omega\left(n^{\frac{1}{2}}\right)$. This finishes the proof of the first part of the theorem. The second part is identical and we thus omit the proof. 
\end{proof}

\section{Heavy and light Hamilton cycles}\label{cubic}
The proof of Corollary \ref{cubic graphs} is based on a simple variant of the HGA and LGA algorithms. Essentially, we will run the same algorithms, merely changing the starting conditions and the definition of $\alpha_i$.
\begin{proof}[Proof of Corollary \ref{cubic graphs}]
Given a graph $G$ and a path $P$ in $G$, we denote by $h(G,P)$ the number of Hamilton \textit{cycles} going through the path $P$. As mentioned before, by \cite{RW92, FJMRW96}, we have that given $G\sim G_{n,3}$, for any $u_0\in V(G)$ \textbf{whp}
\begin{align*}
    \frac{1}{n}\cdot \left(\frac{2}{\sqrt{3}}\right)^n\le h\left(G,\{u_0\}\right) \le n\left(\frac{2}{\sqrt{3}}\right)^n.
\end{align*}
Denoting by $v_1, v_2$ and $v_3$ the neighbours of $u_0$ in $G$, we have by the pigeonhole principle that WLOG $h\left(G,\{u_0v_1\}\right)\ge \frac{2}{3}\cdot h\left(G,\{u_0\}\right)$. Furthermore, given a path $P$ and an end-vertex of the path $u_P$, it has two additional neighbours in $G$, $w_1$ and $w_2$. Since every Hamilton cycle can travel through an edge only once, we then have that $h(G,P)=\alpha\cdot h\left(G,P\cup \{w_1\}\right)+(1-\alpha)h\left(G,P\cup{w_2}\right)$ for some $\alpha\in[0,1]$. We can thus run both the HGA and LGA algorithms in almost the same manner, by initiating the algorithms at their second iteration with $P=\{u_0v_1\}$ and the above definition of $\alpha$. Since at every step we will always choose the edge such that $\alpha>0$, at the end of the process, we obtain a path $P$ of length $n$ that belongs to a Hamilton cycle. 

Noting that $\log_2\left(h(G,\{u_0\})\right)=\left(1+o(1)\right)\left(h(G,\{u_0v_1\})\right)$, we have that the analysis of these variants of the algorithm is precisely the same, and we thus construct a Hamilton cycle with the desired total weight.
\end{proof}

\section{The kernel and long cycles in the giant component of supercritical random graph}\label{discussion}
In his work from 1991 on cycles in the barely supercritical $G(n,p)$, \L{}uczak \cite{L90} showed that the kernel of the giant component is much like a random cubic graph. Ding, Kim, Lubetzky and Peres \cite{DKLP11} and Ding, Lubetzky and Peres \cite{DLP14} studied the anatomy of the giant component of $G(n,p)$, in the barely supercritical regime and the strictly supercritical regime, respectively. Through these results, it can be seen that the study of the giant component of $G(n,p)$ is tightly related to the study of random cubic graphs. When $p=\frac{1+\epsilon}{n}$ and $\epsilon=o(n^{-1/4})$, the description of the giant component is contiguous to the following model. We begin with a random cubic graph on asymptotically $\frac{4\epsilon^3n}{3}$ vertices, which is equivalent to the kernel of the giant. We then assign $i.i.d.$ weights to each edge according to geometric$(\epsilon)$. We subdivide each edge according to the weight assigned, and the resulting graph is the $2$-core of the Giant. We then attach $i.i.d.$ poisson$(1-\epsilon)$, subcritical, Galton-Watson trees to each of the vertices. For larger $\epsilon$, the kernel is contiguous to a multigraph with a random degree distribution, which is dominated by vertices of degree $3$. 

Since there are many standard, fairly good, coupling between geometric$(\lambda)$ and $\frac{1}{\lambda}exp(1)$ random variables, one can see there is a strong connection between our results, and the kernel and the $2$-core of the giant component.

Ajtai, Koml\'os, and Szemer\'edi \cite{AKS81a} proved that a longest cycle in the giant is \textbf{whp} of order $\Theta(\epsilon^2)n$ (see also \cite{KS13} for a simple proof of this fact). The best-known bounds on the length of a longest cycle, as of the writing of this paper, are due to Kemkes and Wormald \cite{KW13}, who showed that a longest cycle has typically length at most $1.739\epsilon^2n$, and the very recent result of Anastos \cite{A22}, who showed that a longest cycle has typically length at least $1.581\epsilon^2n$. Our results imply there is a cycle of asymptotic length at least $\frac{4}{3\epsilon}\cdot 1.1\epsilon^3n=1.46\epsilon^2n$. It is interesting to note that a careful expectation analysis of the number of paths of length $k$, suggests that our algorithm could construct an even longer path if one were to consider paths on only $0.98n$ vertices of the kernel. As one would have to prove log-concentration on the number of paths of such length, and as this would still yield a bound lower than $1.581\epsilon^2n$, we did not venture in that direction. Nevertheless, it is our hope that this algorithm, or insights from it, may help to advance the research towards the sharp asymptotics of a longest cycle in the giant component. It is also worth noting that our explicit (greedy) algorithm process returns a cycle of length typically at most $\frac{4}{3\epsilon}\cdot0.89\epsilon^3n<1.187\epsilon^2n$ covering all the vertices of the kernel of the giant.

Several extensions of the algorithm can be considered. Instead of considering a greedy algorithm, by choosing the best edge at the next immediate step, one could consider a constant number of steps ahead. An advantage in the analysis of this fact is that in the random cubic graph, one would typically expect the $k$-th neighbourhood of a vertex, where $k$ is a constant, to be almost tree-like. On the other hand, there are issues with dependency that arise when one looks beyond the first neighbourhood of the vertex, and the integral calculations become less trivial. Another possible extension would be to greedily construct other structures rather than paths, based on the number of such structures starting from a given vertex. Once again, it is likely that such an extension would lead to a more complicated (and, perhaps, less natural) analysis.

\paragraph{Acknowledgements.} The authors would like to thank Michael Krivelevich for his guidance, advice and many fruitful discussions.
\bibliographystyle{abbrv}
\bibliography{paths}

\begin{thebibliography}{10}

\bibitem{AKS81a}
M.~{Ajtai}, J.~{Koml\'{o}s}, and E.~{Szemer\'{e}di}.
\newblock {The longest path in a random graph}.
\newblock {\em {Combinatorica}}, 1:1--12, 1981.

\bibitem{A01}
D.~Aldous.
\newblock The $\zeta$ (2) limit in the random assignment problem.
\newblock {\em Random Structures Algorithms}, 18(4):381--418, 2001.

\bibitem{A22}
M.~Anastos.
\newblock An improved lower bound on the length of the longest cycle in random
  graphs.
\newblock {\em arXiv preprint arXiv:2208.06851}, 2022.

\bibitem{BFM98}
A.~{Beveridge}, A.~{Frieze}, and C.~{McDiarmid}.
\newblock {Random minimum length spanning trees in regular graphs}.
\newblock {\em {Combinatorica}}, 18(3):311--333, 1998.

\bibitem{DKLP11}
J.~Ding, J.~H. Kim, E.~Lubetzky, and Y.~Peres.
\newblock Anatomy of a young giant component in the random graph.
\newblock {\em Random Structures Algorithms}, 39(2):139--178, 2011.

\bibitem{DLP14}
J.~{Ding}, E.~{Lubetzky}, and Y.~{Peres}.
\newblock {Anatomy of the giant component: the strictly supercritical regime}.
\newblock {\em {European Journal of Combinatorics}}, 35:155--168, 2014.

\bibitem{FJMRW96}
A.~Frieze, M.~Jerrum, M.~Molloy, R.~Robinson, and N.~Wormald.
\newblock Generating and counting hamilton cycles in random regular graphs.
\newblock {\em Journal of Algorithms}, 21(1):176--198, 1996.

\bibitem{FK16}
A.~Frieze and M.~Karo\'{n}ski.
\newblock {\em Introduction to random graphs}.
\newblock Cambridge University Press, Cambridge, 2016.

\bibitem{J99}
S.~Janson.
\newblock One, two and three times log n/n for paths in a complete graph with
  random weights.
\newblock {\em Combinatorics, Probability and Computing}, 8(4):347--361, 1999.

\bibitem{KW13}
G.~Kemkes and N.~Wormald.
\newblock An improved upper bound on the length of the longest cycle of a
  supercritical random graph.
\newblock {\em SIAM Journal on Discrete Mathematics}, 27(1):342--362, 2013.

\bibitem{KS13}
M.~Krivelevich and B.~Sudakov.
\newblock The phase transition in random graphs: {A} simple proof.
\newblock {\em Random Structures Algorithms}, 43(2):131--138, 2013.

\bibitem{L90}
T.~{{\L}uczak}.
\newblock {Component behavior near the critical point of the random graph
  process}.
\newblock {\em {Random Structures Algorithms}}, 1(3):287--310, 1990.

\bibitem{MU17}
M.~Mitzenmacher and E.~Upfal.
\newblock {\em Probability and computing: Randomization and probabilistic
  techniques in algorithms and data analysis}.
\newblock Cambridge University Press, 2017.

\bibitem{RW92}
R.~Robinson and N.~Wormald.
\newblock Almost all cubic graphs are hamiltonian.
\newblock {\em Random Structures and Algorithms}, 3(2):117--125, 1992.

\bibitem{steele1986efron}
J.~M. Steele.
\newblock An {E}fron-{S}tein inequality for nonsymmetric statistics.
\newblock {\em The Annals of Statistics}, 14(2):753--758, 1986.

\bibitem{UP13}
S.~Uryasev and P.~Pardalos.
\newblock {\em Stochastic Optimization: Algorithms and Applications}.
\newblock Applied Optimization. Springer US, 2013.

\end{thebibliography}
\end{document}